\documentclass[fleqn,12 pt]{article}
\usepackage{amsmath,amssymb,amsfonts,amsthm,graphicx}
\usepackage[bookmarksnumbered, colorlinks, plainpages]{hyperref}
\usepackage{cite}
\usepackage{epsfig}
\usepackage{epstopdf}
\usepackage{footnote}
\theoremstyle{plain}
\newtheorem{theorem}{Theorem}[section]
\newtheorem{lemma}[theorem]{Lemma}

\theoremstyle{definition}
\newtheorem{definition}[theorem]{Definition}

\date{}

\title{On the distance from a normal matrix polynomial to the set of matrix polynomials with a prescribed multiple eigenvalue}

\author{E. Kokabifar $^1$, G.B. Loghmani $^2$\\
\footnotesize{$^{1,2}$Faculty of Mathematics, Yazd University, Yazd, Iran.}\\
 \footnotesize{e-mail:$^1$ e.kokabifar@yahoo.com, $^2$ loghmani.yazd.ac.ir.}\\
}
\date{}

\begin{document}

\maketitle
%\subjclass[2010]{Primary 15A18; Secondary 65F35.}
%\keywords{Matrix; Singular value; Eigenvalue; Perturbation}\\
%\indent $^{*}$ Corresponding auther}

\begin{abstract}
This paper concerns the bounds for spectral norm distance from a normal matrix polynomial $P(\lambda)$ to the set of matrix polynomials that have $\mu$ as a multiple eigenvalue. Also construction of associated perturbations of $P(\lambda)$ is considered.
\end{abstract} \maketitle

\section{Introduction}
For a complex $n \times n$ matrix $A$ and a given complex number $\mu$, the spectral norm distance from $A$ to the matrices that have $\mu$ as a multiple eigenvalue was proved by A. M. Malyshev \cite{malyshev}. Malyshev's results extended by Lippert \cite{lipert} and Gracia \cite{gracia} and they compute a 2-norm distance from $A$ to the set of matrices with two  prescribed eigenvalue. In 2004, Ikramov and Nazari \cite{ikramovasli} show that Malyshev's formula not a viable method for the case of normal matrices. Also an operative manner for this case was introduced by them. Moreover, Nazari and Rajabi \cite{nazarirajabi} noted this issue for distance from a normal matrix $A$ to the set of matrices with two prescribed eigenvalues.

In 2008, a spectral norm distance from a matrix polynomial $P(\lambda)$ to the matrix polynomials that have $\mu$ as an eigenvalue of geometric multiplicity at least $\kappa$, and a distance from $P(\lambda)$ to the matrix polynomials that have $\mu$ as a multiple eigenvalue was introduced by Papathanasiou and Psarrakos \cite{papa} . They computed the first distance and also obtained bounds for the second one, constructing an associated perturbations of $P(\lambda)$. In this paper, at first we illustrate by an example that procedure described in \cite{papa} not a efficient method for the case of normal matrix polynomials. Then a suitable process is presented.
\section{Preliminaries}
In this section some definitions of matrix polynomials are reviewed. A good reference for the theory of matrix polynomials is \cite{lancaster}.  Moreover, we illustrate by a specific example that method described in \cite{papa} is useless for normal matrix polynomials.
\begin{definition}
For $A_j\in\mathbb{C}^{n \times n}(j = 0,1,...,m)$ with det$(A_m)\neq 0$ and a complex variable $\lambda$, we define the matrix polynomial $P(\lambda)$ as
\begin{equation}\label{plambda}
P(\lambda ) = A_m \lambda ^m  + A_{m - 1} \lambda ^{m - 1}  + ... + A_1 \lambda  + A_0.
\end{equation}
\end{definition}
If for a scalar $\mu$ and some nonzero vector $\upsilon \in {\mathbb{C}^{n}}$, it holds that $P(\mu)\upsilon=0,$ then the scalar $\mu$ called an \textit{eigenvalue} of $P(\lambda)$ and the vector $\upsilon$ is known as a \textit{right eigenvector} of $P(\lambda)$ corresponding to $\mu$. Similarly, nonzero vector $\nu  \in {\mathbb{C}^{n}}$ is known as a left eigenvector of $P(\lambda)$ corresponding to $\mu$ if we have $\nu^*P(\mu)=0$. Multiplicity of $\mu$ as a root of the scalar polynomial det$P(\lambda)$ called as its \textit{algebraic multiplicity} and number of linear independent eigenvectors corresponding to $\mu$ is \textit{geometric multiplicity}. algebraic multiplicity of an eigenvalue is always greater or equal to its geometric multiplicity. If algebraic and geometric multiplicities of an eigenvalue are equal, this eigenvalue called \textit{semisimple}, otherwise it named \textit{defective}.
\begin{definition}
Let $P(\lambda)$ be a matrix polynomial as in (\ref{plambda}), if there exists a unitary matrix $U \in \mathbb{C}^{n \times n}$ such that $U^*P(\lambda)U$ is diagonal for all $\lambda \in \mathbb{C}$, then $P(\lambda)$ named \textit{weakly normal}. Moreover, $P(\lambda)$ is called \textit{normal} if all the eigenvalues of $P(\lambda)$ are semisimple.
\end{definition}
For the matrix polynomials of the form $A+\lambda B$, two concepts weakly normal and normal are coincide. The matrix polynomial $P(\lambda)$ is weakly normal if and only if for every $\mu\in \mathbb{C}$ the matrix $P(\mu)$ is a normal matrix. Moreover, all of the coefficient matrices $A_i, (i=1,\ldots, m)$ in (\ref{plambda}) also all linear combinations of them are normal matrices \cite{papanormal}.

Here, some of results obtained in sections $4$ and $5$ of \cite{papa} are reviewed briefly.

For the matrix polynomial $P(\lambda)$ in (\ref{plambda}) and a complex variable $\gamma$, Papathanasiou and Psarrakos \cite{papa} define
\[F\left[ {P(\lambda );\gamma } \right] = {\left[ {\begin{array}{*{20}{c}}
{P(\lambda )}&0\\
{\gamma P'(\lambda )}&{P(\lambda )}
\end{array}} \right]_{2n \times 2n}},\]
where $P'(\lambda )$ denotes the derivative of $P(\lambda )$ with respect to $\lambda$.
\begin{lemma}\textsc{\cite{papa}}.\label{dotaee}
Let $\mu \in \mathbb{C}$ and $\gamma_*> 0$ be a point where the singular value $s_{2n -
1} (F[P(\mu );\gamma ])$ attains its maximum value, and $s_*=s_{2n -
1} (F[P(\mu);\gamma_* ])>0$. Then there exist a pair $\left[
{\begin{array}{*{20}c}
   {u_1 (\gamma _* )}  \\
   {u_2 (\gamma _* )}  \\
\end{array}} \right], \left[ {\begin{array}{*{20}c}
   {v_1 (\gamma _* )}  \\
   {v_2 (\gamma _* )}  \\
\end{array}} \right] \in\mathbb{C} ^{2n}~( u_k (\gamma _* ),v_k (\gamma _* ) \in \mathbb{C}^n,~k =
1,2)$ of left and right singular vectors of $s_*$ respectively,
such that

$1.~u^* _2 (\gamma _* )P'(\mu)v_1 (\gamma _* ) =0, $ and

$2.$~the $n\times 2$ matrices $ U(\gamma _* ) = [u_1 (\gamma _*
)~u_2 (\gamma _* )]_{n \times 2}$ and $ V(\gamma _* ) = [v_1
(\gamma _* )~v_2 (\gamma _* )]_{n \times 2} $ satisfy $
U^*(\gamma _* ) U(\gamma _* ) = V^* (\gamma _* ) V(\gamma _* ).$
\end{lemma}
Suppose that weights $w = \{ \omega _0 ,\omega _1 ,...,\omega _m \}$ are given, such that $w$ is a set of nonnegative coefficients with $\omega _0 > 0$. The scalar polynomial $w(\lambda )$ corresponding to the weights is defined in the form
\begin{equation*}
w(\lambda ) = {w_m}{\lambda ^m} +  \cdots  + {w_1}\lambda  + {w_0}.
\end{equation*}
Moreover, consider the matrix
\begin{equation*}
{\Delta _{{\gamma _*}}} =  - {s_*}U({\gamma _*})\left[ {\begin{array}{*{20}{c}}
1&{ - {\gamma _*}\phi }\\
0&1
\end{array}} \right]V{({\gamma _*})^\dag },
\end{equation*}
where ${V}({\gamma_*})^\dag$ is the
\emph{Moore-Penrose pseudoinverse} of ${V}({\gamma_*} )$ and the quantity $\phi$ is
\begin{equation*}
\phi  = \frac{{w'(\left| \mu  \right|)}}{{w(\left| \mu  \right|)}}\frac{{\bar \mu }}{{\left| \mu  \right|}},
\end{equation*}
if $\mu=0$, then by convention we set $\frac{{\bar \mu }}{{\left| \mu  \right|}}=0.$

If $\gamma_* >0$ and $\mu  \in {\mathbb{C}}$ is not an eigenvalue of $P'(\lambda)$, then \cite[Theorem 19]{papa} implies that matrix polynomial $Q_{\gamma_*}(\lambda)=P(\lambda)+\Delta_{\gamma_*}(\lambda)$ has $\mu$ as a defective eigenvalue. Where
\begin{equation*}
{\Delta _{{\gamma _*}}}(\lambda ) = \sum\limits_{j = 0}^m {\left( {\frac{{{w_j}}}{{w(\left| \mu  \right|)}}{{\left( {\frac{{\bar \mu }}{{\left| \mu  \right|}}} \right)}^j}{\Delta _{{\gamma _*}}}} \right){\lambda ^j}}.
\end{equation*}
Now for a specific example, let us consider the normal matrix polynomial $P(\lambda)$ as mentioned in \cite[section 3]{papanormal} of the form
\begin{equation}\label{pnormal}
P(\lambda ) = \left[ {\begin{array}{*{20}{c}}
1&0&0\\
0&1&0\\
0&0&1
\end{array}} \right]{\lambda ^2} + \left[ {\begin{array}{*{20}{c}}
{ - 3}&0&0\\
0&{ - 1}&0\\
0&0&3
\end{array}} \right]\lambda  + \left[ {\begin{array}{*{20}{c}}
2&0&0\\
0&0&0\\
0&0&2
\end{array}} \right].
\end{equation}
Consider the set of weights $w = \left\{ {1,1,1} \right\}$ and $\mu=3$. By applying the procedure described sections $4$ and $5$ of \cite{papa} (and reviewed briefly above), it can obtained that $s_{5}(F[P(3);\gamma])$ attains its maximum value at $\gamma_{*}=1$, and $s_*=s_{5}(F[P(3);\gamma_*])=4$.
Also we have find the $Q(\gamma_*)$ as following
\begin{equation*}
{Q_{{\gamma _*}}}(\lambda ) = \left[ {\begin{array}{*{20}{c}}
1&0&0\\
0&{0.6013}&0\\
0&0&1
\end{array}} \right]{\lambda ^2} + \left[ {\begin{array}{*{20}{c}}
{ - 3}&0&0\\
0&{ - 1.3987}&0\\
0&0&3
\end{array}} \right]\lambda  + \left[ {\begin{array}{*{20}{c}}
2&0&0\\
0&{ - 0.3987}&0\\
0&0&2
\end{array}} \right],
\end{equation*}
it is straightforward to see that $\mu=3$ is not an eigenvalue of ${Q_{{\gamma _*}}}(\lambda )$.
% are listed below
%\[\sigma \left( {{Q_{{\gamma _*}}}(\lambda )} \right) = \left\{ {1,2, - 2, - 1,2.5830, - 0.2567} \right\},\]
Furthermore,
\begin{equation*}
u_2^*({\gamma _*})P'(\mu ){v_1}({\gamma _*}) = {\rm{  - 1}}{\rm{.5385}},\hspace{.75cm} {\rm and} \hspace{.75cm} {\left\| {{U^*}({\gamma _*})U({\gamma _*})-{V^*}({\gamma _*})V({\gamma _*})} \right\|_2} = 0.3846,
\end{equation*}
imply that none of propositions of the Lemma \ref{dotaee} is not confirmed.
\section{Normal matrix polynomial}
As was mentioned in the previous section, described method in \cite{papa} is not efficient for the case of normal matrix polynomials. Moreover, both results of Lemma \ref{dotaee} are violated.  In fact violation of the first part of Lemma \ref{dotaee} is cause of violation of the second part.
By similar analysis fulfilled in \cite{onaremarkable} it can be showed that when $\gamma \geq 0$ started to rise $s_{2n-1} (F[P(\mu);\gamma ])$ increases and $s_{2n-2} (F[P(\mu);\gamma ])$ decreases, and where to next this process will reversed. Consequently, there exists a point such as $\gamma_*$ where $s_{2n-1} (F[P(\mu);\gamma ])$ attains its maximum value. Furthermore, at $\gamma=\gamma_*$ we have
$s_*=s_{2n-1} (F[P(\mu);\gamma_* ])=s_{2n-2} (F[P(\mu);\gamma_* ])$, that means $s_*$ is a multiple singular value of $F[P(\mu);\gamma ]$. The graphs of the $s_{5}(F[P(3);\gamma])$  and $s_{4}(F[P(3);\gamma])$ for $\gamma\in [0, 10]$ are plotted in Fig 1. 
\begin{figure}
\centering
\includegraphics[width=0.52\linewidth]{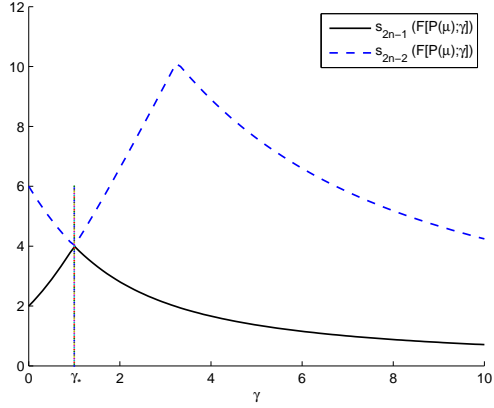}
\caption{$s_{2n-1} (F[P(\mu);\gamma ])$ and $s_{2n -2} (F[P(\mu);\gamma])$}
\label{fig:fig1}
\end{figure}

Suppose that ${u^{2n - 1}}({\gamma _*}), {v^{2n - 1}}({\gamma _*})$ and ${u^{2n - 2}}({\gamma _*}), {v^{2n - 2}}({\gamma _*})$ are a pair of singular vectors of $s_{2n-1} (F[P(\mu);\gamma_* ])$ and $s_{2n-2} (F[P(\mu);\gamma_* ])$, respectively. As we have seen, computation of first proposition in the Lemma \ref{dotaee} for $s_{2n-1} (F[P(\mu);\gamma_* ])$ yields
\begin{equation}\label{manfi}
{u_2}^{(2n - 1)}{({\gamma _*})^*}P'(\mu ){v_1}^{(2n - 1)}({\gamma _*}) = {\rm{  - 1}}{\rm{.5385}},
\end{equation}
also similar computations but for $s_{2n-2} (F[P(\mu);\gamma_* ])$ leads to
\begin{equation}\label{mosbat}
{u_2}^{(2n - 2)}{({\gamma _*})^*}P'(\mu ){v_1}^{(2n - 2)}({\gamma _*}) = 2.4,
\end{equation}
According to Lemma 16 of \cite{papa}, it is straightforward to see that ${u_2}{({\gamma })^*}P'(\mu ){v_1}({\gamma })$ can be explained as the derivative of the corresponding singular value with respect to $\gamma$. Therefore, negativity and positivity of the numbers in (\ref{manfi}) and (\ref{mosbat})(respectively) means that $s_{2n-1} (F[P(\mu);\gamma ])$ and $s_{2n-2} (F[P(\mu);\gamma ])$ are decreasing and increasing functions, respectively.

Hereafter, we are looking for a pair of vectors $u$ and $v$ in the form
\begin{equation}\label{uvnew}
u({\gamma _*}) = \alpha {u}^{(2n - 1)}({\gamma _*}) + \beta {u}^{(2n - 2)}({\gamma _*}),\hspace{.59cm} {\rm and}\hspace{.59cm} v({\gamma _*}) = \alpha {v}^{(2n - 1)}({\gamma _*}) + \beta {v}^{(2n - 2)}({\gamma _*}),
\end{equation}
such as satisfy 
\begin{equation}\label{upvnew}
{u_2}{({\gamma _*})^*}P'(\mu ){v_1}({\gamma _*}) =0.
\end{equation}
Where for the two scalars $\alpha$ and $\beta$ we have ${\left| \alpha  \right|^2} + {\left| \beta  \right|^2} = 1$.
For doing this, substituting $u$ and $v$ in (\ref{uvnew}) into (\ref{upvnew}) leads to
\begin{equation}\label{alphabeta}
\left[ {\begin{array}{*{20}{c}}
{\bar \alpha }&{\bar \beta }
\end{array}} \right]M\left[ {\begin{array}{*{20}{c}}
\alpha \\
\beta 
\end{array}} \right] = 0,
\end{equation}
where
\begin{equation*}
M = \left[ {\begin{array}{*{20}{c}}
{{u_2}^{(2n - 1)}{{({\gamma _*})}^*}P'(\mu ){v_1}^{(2n - 1)}({\gamma _*})}&{{u_2}^{(2n - 1)}{{({\gamma _*})}^*}P'(\mu ){v_1}^{(2n - 2)}({\gamma _*})}\\
{{u_2}^{(2n - 2)}{{({\gamma _*})}^*}P'(\mu ){v_1}^{(2n - 1)}({\gamma _*})}&{{u_2}^{(2n - 2)}{{({\gamma _*})}^*}P'(\mu ){v_1}^{(2n - 2)}({\gamma _*})}
\end{array}} \right].
\end{equation*}
It easy to see that $M$ is a indefinite Hermitian matrix, which implies that there exists a nontrivial solution for (\ref{alphabeta}). Suppose that the matrix $M$ has a spectral decomposition of the form
\begin{equation*}
M = U\left[ {\begin{array}{*{20}{c}}
{{\lambda _1}}&0\\
0&{{\lambda _2}}
\end{array}} \right]{U^*},
\end{equation*}
assume a unit vector $\left[ {\begin{array}{*{20}{c}}
\xi \\
\eta 
\end{array}} \right]$ and set $\left[ {\begin{array}{*{20}{c}}
\alpha \\
\beta 
\end{array}} \right] = U\left[ {\begin{array}{*{20}{c}}
\xi \\
\eta 
\end{array}} \right]$. So, the equation (\ref{alphabeta}) turns on 
\[\left[ {\begin{array}{*{20}{c}}
{\bar \xi }&{\bar \eta }
\end{array}} \right]\left[ {\begin{array}{*{20}{c}}
{{\lambda _1}}&0\\
0&{{\lambda _2}}
\end{array}} \right]\left[ {\begin{array}{*{20}{c}}
\xi \\
\eta 
\end{array}} \right] = 0.\]
Finally $\left[ {\begin{array}{*{20}{c}}
\xi \\
\eta 
\end{array}} \right]$ and consequently $\left[ {\begin{array}{*{20}{c}}
\alpha \\
\beta 
\end{array}} \right]$ can be obtain form
\[{\left| \xi  \right|^2}{\lambda _1} + {\left| \eta  \right|^2}{\lambda _2} = 0,\qquad {\rm and}\qquad{\left| \xi  \right|^2} + {\left| \eta  \right|^2} = 1,\]
that straightforward yields
\[\xi  = \sqrt {\frac{{\left| {{\lambda _1}} \right|}}{{\left| {{\lambda _1}} \right| + \left| {{\lambda _2}} \right|}}} ,\qquad{\rm and}\qquad \eta  = \sqrt {\frac{{\left| {{\lambda _2}} \right|}}{{\left| {{\lambda _1}} \right| + \left| {{\lambda _2}} \right|}}}. \]
Now, return to the above example of normal matrix polynomial. By applying what is discussed for the normal matrix polynomial $P(\lambda)$ in (\ref{pnormal}) we have
\begin{equation*}
\alpha  = {\rm{ - 0}}{\rm{.6250}},\qquad{\rm and}\qquad\beta  = {\rm{ - 0}}{\rm{.7806}},
\end{equation*}
also the two vectors $u({\gamma _*})$ and $v({\gamma _*})$ in (\ref{uvnew}) satisfy
\[u_2^*({\gamma _*})P'(\mu ){v_1}({\gamma _*}) = {\rm{  - 2}}{\rm{.2204}} \times {10^{ - 16}},\]
and
\[{\left\| {{U^*}({\gamma _*})U({\gamma _*}) - {V^*}({\gamma _*})V({\gamma _*})} \right\|_2} = {\rm{3}}{\rm{.3479}} \times {10^{ - 16}},\]
Moreover, the matrix polynomial $Q_{\gamma_*}(\lambda)$ that has $\mu=3$ as a multiple eigenvalue can be find as following
\begin{eqnarray*}
{Q_{{\gamma _*}}}(\lambda ) = \left[ {\begin{array}{*{20}{c}}
{0.7722}&{ - 0.0955}&0\\
{0.0527}&{0.6065}&0\\
0&0&1
\end{array}} \right]{\lambda ^2} &+& \left[ {\begin{array}{*{20}{c}}
{ - 3.2278}&{ - 0.0955}&0\\
{0.0527}&{ - 1.3935}&0\\
0&0&3
\end{array}} \right]\lambda \\ &+& \left[ {\begin{array}{*{20}{c}}
{1.7722}&{ - 0.0955}&0\\
{0.0527}&{ - 0.3935}&0\\
0&0&2
\end{array}} \right].
\end{eqnarray*}

\end{document}